\documentclass[a4paper]{article}
\usepackage{graphicx}
\usepackage{epstopdf}
\usepackage{amsfonts}
\usepackage{amsmath}

\begin{document}

\title{Accurate Spectral Collocation Solutions to some Bratu's Type Boundary Value Problems.}

\author{C\u{a}lin-Ioan Gheorghiu
\thanks{Tiberiu Popoviciu Institute of Numerical Analysis, Str. Fantanele 57, Cluj-Napoca, Romania,
email: ghcalin@ictp.acad.ro
}}

\date{ }
\maketitle

\begin{abstract}
We solve by Chebyshev spectral collocation some genuinely nonlinear Liouville-Bratu-Gelfand type, 1D and a 2D boundary value problems. The problems are formulated on the square domain $[-1, 1]\times[-1, 1]$ and the boundary condition attached is a homogeneous Dirichlet one. We pay a particular attention to the bifurcation branch on which a solution is searched and try to estimate empirically the attraction basin for each bifurcation variety. 

The first eigenvector approximating the corresponding the first eigenfunction of the linear problem is used as an initial guess in solving the nonlinear algebraic system of Chebyshev collocation to find the \textquotedblleft{small}\textquotedblright solution. For the same value of the bifurcation parameter we use another initial guess, namely lowest basis function (1 point approximation), to find the \textquotedblleft{big}\textquotedblright solution. The Newton-Kantorovich method solves very fast the nonlinear algebraic system in no more than eight iterations. Beyond being exact, the method is numerically stable, robust and easy to implement. Actually, the MATLAB code essentially contains three programming lines. It by far surpasses in simplicity and accuracy various methods used to solve some well-known problems. 

We end up by providing some
numerical and graphical outcomes in order to underline the validity and the effectiveness of our method, i.e., norms of Newton updates in solving the algebraic systems and the decreasing rate of Chebyshev coefficients of solution.


\end{abstract}

\section{Introduction} \label{intro}
The boundary value problem attached to the Liouville-Bratu-Gelfand operators in two dimensions is classical in mathematics and physics. Nevertheless, computational methods for highly accurate solutions are still of much current interest, particularly in applications.
Thus for technical point of view, the exothermic chemical reaction of bulk materials held in a storage is a common hazard. 

Consider the distribution of temperature, $u,$ within such a material, which is increased by the release of energy from reaction. We have the model (see for instance \cite{Grindrod})
\begin{equation}
u_{t}=D\Delta u+\alpha e^{u},\ \ x\in \Omega ,\ \ t>0, \label{exothermic}
\end{equation}
where $D$ is the heat conductivity and$\alpha$ is a positive constant. The rate $\alpha e^{u}$ is actually a simplification of the more usual Arrhenius law.
Let us suppose that $u$ is held fixed, say $=u_{b},$ at the boundary of the domain $\Omega.$ This is the most representative example in the class of problems we are going to consider.

Various numerical methods, finite elements, finite difference and finite volumes for space discretization along with various finite difference schemes to march in time, have been designed in order to solve initial-boundary value problems attached to the equation (\ref{exothermic}). We will mainly confine ourselves to one and two dimensional steady equations coming from (\ref{exothermic}). We will also solve some nonlinear problems where the non linearity in (\ref{exothermic}) is of the form $f(\varepsilon,u):=f(\varepsilon, e^{u})$ with $f:{\mathbb R}_{+} \times {\mathbb R} \rightarrow {\mathbb R}$ continuous and differentiable with respect to both arguments. 

Our main aim is to show that the Chebyshev spectral collocation (ChC) along with Newton-Kantorovich produces highly accurate solutions to these problems irrespective to the non linearity and corner singularity of the domain.
The method is robust, very efficient and easy implementable. For a given value of bifurcation parameter we compute both solutions, the \textquotedblleft{small}\textquotedblright and the \textquotedblleft{big}\textquotedblright one. For 2D case, up to our knowledge this a novelty.

The strategy in question emerges from the work of John Boyd \cite{Boyd}. However, recently using ChC along with Newton-Kantorovich we have successfully solved a third-order singular and nonlinear boundary value problem on the half-line (see \cite{Gheorghiu}). Actually, we have used o Chebfun code in order to implement both algorithms. This problem possesses an integral invariant and thus we could guarantee an accuracy tending to the machine precision. For some similar nonlinear boundary value problems we have used, beyond collocation method, Galerkin type methods in order to solve them. Thus we quote some of them \cite{BG}, \cite{GT1} and \cite{GT}.

In a recent paper \cite{Kouibia} the authors propose a variational method in order to solve Bratu’s problem for
two dimensions in an adequate approximate space of bi quadratic spline functions. Their numerical results look correct but the authors do not provide very convincing arguments in support of their accuracy.  

This paper divides naturally into two parts. The first Section \ref{1DB} is devoted to the 1D Bratu's problem. We find the \textquotedblleft{small}\textquotedblright and the \textquotedblleft{big}\textquotedblright solution and study their stability. At least graphically we show that the Newton-Kantorovich process is of order larger than two. We also observe that for both solutions their Chebyshev coefficients decrease exponentially. A simple Chebfun code has been very efficient. In the second Section \ref{2DB}, in principle, we have done almost the same for the 2D Bratu's problem. However, in order to discretize the 2D Bratu's equation we have used the classical Chebyshev differentiation matrices. Moreover, we believe that we are actually showing for the first time how the Chebyshev coefficients of solutions to boundary or eigenvalues problems decrease.

Some concluding remarks are also supplied.

\section{1D Bratu's problem} \label{1DB}
Suppose now $\Omega:=[-L, L]$ in one dimension and consider the steady-state problem
\begin{equation}
u_{xx}+\lambda e^{u}=0,\ u\left( -L\right) =u\left( L\right) =0,\ L>0, \label{1DB_eq}
\end{equation}
where $\lambda :=\alpha /D$. The solution of this problem is even, so normally one solves the half-problem on $[0, 1]$ under mixed boundary conditions $u^{\prime }\left( 0\right) =u\left( L\right) =0.$ Gelfand in \cite{Gelfand} solved Bratu's problem as an initial value problem with conditions $u^{\prime }\left( 0\right) =0,$ $u\left( 0\right) =A,$ $A>$ by shooting. He found a solution of the form%
\begin{equation}
w\left( x\right) =2\ln \left[ \frac{e^{A/2}}{\cosh \left( Bx\right) }\right]
,\ B:=\left( \lambda e^{A}/2\right) .\label{w_sol}
\end{equation}
For a fixed $L,$ by enforcing the boundary condition at $x:=L$ in solution (\ref{w_sol}) one can find the following transcendental relation between $A$ and $\lambda$ 
\begin{equation}
\frac{e^{A/2}}{cosh\left( \sqrt{\lambda e^{A}/2}L\right)} =1.  \label{bifurc_eq}
\end{equation}
It provides the bifurcation pattern for problem (\ref{1DB_eq}) (see again Gelfand \cite{Gelfand} for other comments on global solutions to Bratu's problem (\ref{1DB_eq}) and their invariant). 

\subsection{Bifurcation of positive solutions} \label{1DB_bifurc}
\begin{figure}
  \includegraphics{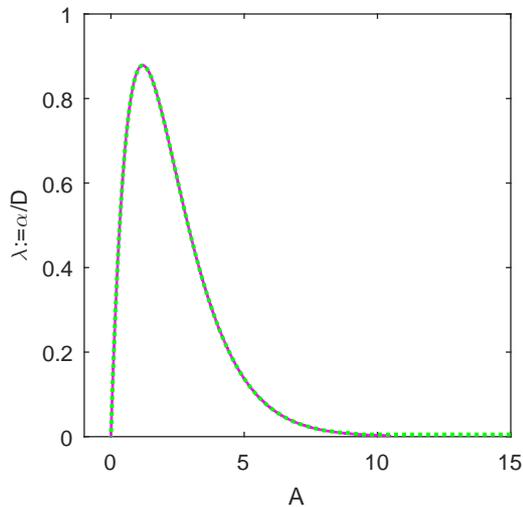}
\caption{Bifurcation diagram for 1D Bratu's problem.}
\label{fig:1}       
\end{figure}
In order to visualize the bifurcation pattern (\ref{bifurc_eq}) we use the \texttt{roots} method for class Chebfun2 (see \cite{DHT} Part II Functions of two variables (Chebfun2)). Thus we have got Fig. \ref{fig:1} with the critical point (which is defined to be where $\lambda(A)$ is a maximum) of coordinates $$\left(A_{*}, \lambda_{*}\right)=\left(1.187331536443172, 0.8786312538512331\right).$$ Specifically, solutions are unique for $\lambda \leq 0$ and for a single positive value $\lambda_{*}$ called the Frank-Kamenetskii critical value; do not exist for $\lambda >\lambda_{*}$; and two solutions exist for  $0<\lambda <\lambda_{*}.$ For $\lambda_{*}$ Boyd reported in \cite{Boyd1} the value $0.8784576797812903015$ which is fairly closed to that reported above.
We have to mention that for $L:=1/2$ the critical value $\lambda_{*}=3.51360308$ and for $L:=2$ we have got $\lambda_{*}=0.2196644$ which means this value decreases as the length of the integration interval increases.

\subsection{Chebfun solutions to 1D Bratu's problem and their linear stability.} \label{1DB-solution}
 \begin{figure*}
  \includegraphics[width=0.75\textwidth]{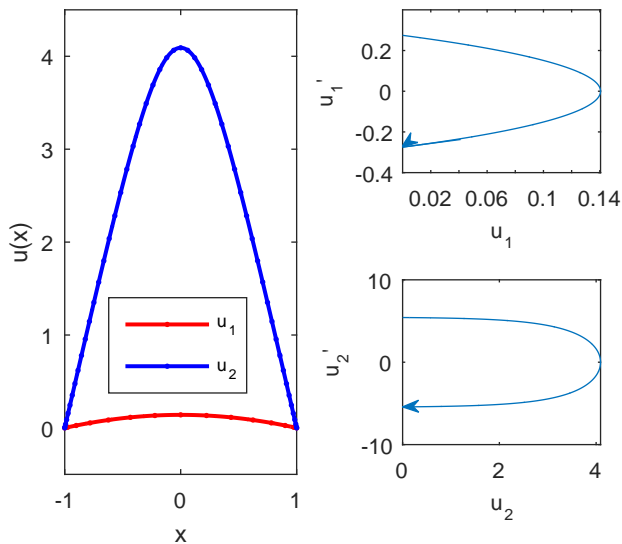}
\caption{Solutions $u_{1}$ and $u_{2}$ (left panel) and their representation in the phase plane (right panels).}
\label{fig:2}       
\end{figure*}
For the same value of $\lambda$, i.e., $\lambda:=0.25$ we compute, using a simple Chebfun code, the solutions corresponding to the two branches of bifurcation diagram. 
The first one, denoted by $u_{1},$ results from Chebfun’s default
initial guess (the zero function). This is the \textquotedblleft{small}\textquotedblright solution. The other, the \textquotedblleft{big}\textquotedblright one denoted by $u_{2},$ results from the alternative initial guess
$6(1-x^2)$, i.e., the simplest basis function for Galerkin spectral method. They both are represented in the left panel of Fig. \ref{fig:2} along with their trajectories in the phase plane (right panels).
\begin{figure*}
  \includegraphics[width=0.75\textwidth]{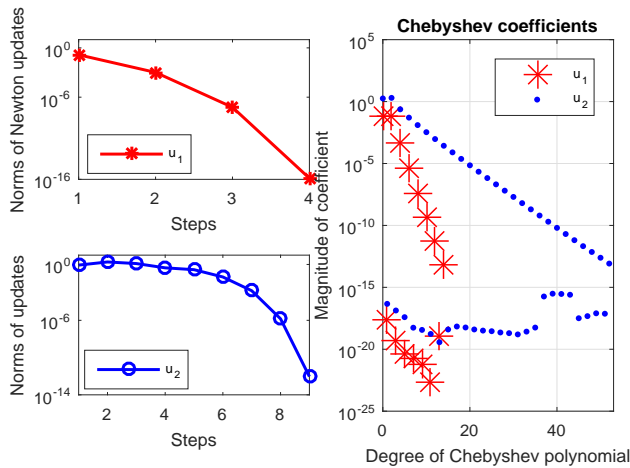}
\caption{Norms of Newton updates in solving for $u_{1}$ and $u_{2}$ (left panel) and the Chebyshev coefficients of these solutions (right panel).}
\label{fig:3}       
\end{figure*}

In order to evaluate the speed of the convergence process we display in the left panels of Fig. \ref{fig:3} the evolution of Newton's updates in computing both solutions. In the right panel, in a log-linear plot, we display the decreasing rate of the even Chebyshev coefficients of solutions. The odd coefficients are lower than $O(10^{-15}).$ A brief look to Fig. \ref{fig:3} leads us to the conclusion that the computation process performed by Chebfun code is much more expensive for the second solution than for the first one (the number of Newton's iterations and the number of Chebyshev coefficients of solutions) is considerably higher for $u_{2}$ than for $u_{1}$.

In addition, we will show in what follows that the first \textquotedblleft{small}\textquotedblright solution is stable but the second one \textquotedblleft{big}\textquotedblright is not. 
Let's consider now the unsteady problem corresponding to (\ref{1DB_eq}).
It reads%
\begin{equation}
u_{t}=u_{xx}+\lambda e^{u},\ x\in \left( -L,\ L\right) ,\ u\left( \pm
L\right) =0,  \label{1DB_eq_unsteady}
\end{equation}%
where $\lambda $ is now fixed in the interval $(0,\ \lambda _{\ast }).$ 

Linearizing about $u=u_{i},\ i=1,2$ we have%
\begin{equation*}
v_{t}=v_{xx}+\left( \lambda e^{u_{i}}\right) v,\ x\in \left( -L,\ L\right)
,\ u\left( \pm L\right) =0,\ i=1,2.
\end{equation*}%
Now the eigenvalue problems%
\begin{equation*}
-\phi _{xx}-\left( \lambda e^{u_{i}}\right) \phi =\mu \phi ,\ x\in \left(
-L,\ L\right) ,\ \phi \left( \pm L\right) =0,\ i=1,2,
\end{equation*}%
have the eigenvalues%
\begin{equation*}
\mu _{i,n}=\left( \frac{n\pi }{2L}\right) ^{2}-\lambda e^{u_{i}},\ i=1,2.
\end{equation*}%
It easy to verify that for $L:=1,\ \mu _{1,n}>0$ and thus solution $u_{1}$is
stable, but $\mu _{2,n}<0$ which means that solution $u_{2}$ is unstable. 

It is important to observe that in the text \cite{Trefethen} the authors obtain the bifurcation diagram (\ref{bifurc_eq}) by integrating successively the problem (\ref{1DB_eq}) using a Chebfun code.

 \section{2D Bratu's problem} \label{2DB}
 
The two-dimensional Bratu’s problem is an elliptic partial differential equation
with homogeneous Dirichlet boundary conditions. This problem is given by
\begin{equation}
\left\{ 
\begin{array}{c}
\Delta u+\lambda e^{u}=0,\ \textrm{in}\ \Omega :=\left( -L,\ L\right) \times\left( -L,\
L\right) ,\ L>0, \\ 
u=0,\ \textrm{on}\ \partial \Omega .%
\end{array}%
\right.   \label{2DB_eq}
\end{equation}
The problem is a
nonlinear eigenvalue one that is commonly used as a test problem for many numerical methods.

\subsection{Bifurcation of positive solutions} \label{2DB_bifurc}

In his paper \cite{Boyd} John Boyd finds in a very ingenious and at the same time a simple way an approximate  relationship between the eigenparameter $\lambda$ and the maximum norm of solution denoted by $A$. He uses the lowest basis function that is a polynomial and vanishes at the four sides of the square, namely 
\begin{equation}
u\left( x,y\right) :=A\left( 1-x^{2}\right) \left( 1-y^{2}\right) ,\ A\in \mathbb{R}. \label{1point_approx}
\end{equation}
Then he uses a type of the weighted residual argument and finds the one-point analytic approximation 
\begin{equation}
\lambda \approx 3.2Ae^{-0.64A}.  \label{1point_eigenvalue_approx}
\end{equation}
of the bifurcation diagram.
He improves this relationship working with a three-point collocation approximation. 
In any case, the bifurcation diagram remains perfectly analogous to that in Fig. \ref{1DB_bifurc}.

\paragraph{Corner singularities and symmetries.} \label{sing_sim}

In the same paper J. Boyd also observe that the solution to problem (\ref{2DB_eq}) has singularities at the four corners of
the domain proportional to%
\begin{equation*}
r^{2}\log \left( r\right), 
\end{equation*}
where $r$ is the radial coordinate of a local polar coordinate system
centered on a given corner. This is a very weak singularity in the sense
that second-order finite differences still give an error that is
asymptotically of order $h^{2}$, where $h$ is the grid spacing, as noted in
Haidvogel and Zang \cite{HZ}.

It is clear that the solution to Bratu's equation is symmetric about the origin in both
$x$ and $y$. Consequently, \textit{only the even-degree Chebyshev polynomials are
needed to construct the appropriate two-dimensional basis functions}. For a
given resolution, this double parity reduces the size of the basis set by a
factor of four. Anyway, we will not exploit this symmetry.
However, if we rotate the square through an angle of $90^{o}$ neither the equation
nor the boundary conditions are altered. Thus, in the language of group theory, the solution is also 
\textquotedblleft{invariant under the rotation group $C_{4}$}\textquotedblright
which is the group of rotations through any integer multiple of $90^{o}$.
 
\subsection{ChC solutions to 2D Bratu's problem.}
\label{2DB-solution}
In order to find this solution we have used the ChC method based on the collocation differentiation matrices from \cite{WR}.

This procedure starts with the second order Chebyshev differentiation matrix $\mathbf{D%
}^{\left( 2\right) }$ of order $N.$ In order to implement the homogeneous
Dirichlet boundary conditions we discard the first and the last rows and
columns of this matrix and obtain the square matrix of order $N-2$ denoted by $%
\widetilde{\mathbf{D}}^{\left( 2\right).}$ Now in order to obtain the
differentiation matrices corresponding to the second order partial
derivatives $\partial ^{2}u/\partial x^{2}$ and $\partial ^{2}u/\partial
y^{2}$ we take the following Kronecker products 
\begin{equation*}
\mathbf{D20}=kron(\mathbf{eye}(N-2,N-2),\widetilde{\mathbf{D}}^{\left(
2\right) });\mathbf{D02}=kron(\widetilde{\mathbf{D}}^{\left( 2\right) },%
\mathbf{eye}(N-2,N-2)).
\end{equation*}%
The discrete Laplacian is now simply $\mathbf{\Delta =D20+D02}$, and the
MATLAB code \texttt{eigs} is now used to find the set of eigenvalues and
eigenfunctions. Actually using the code \texttt{reshape} we transform eigenvector 
$U$ of length $N^2$ into the eigensolution matrix $\mathbf{UU}$ of dimension $N \times N.$

When we solve the 2D Bratu's problem we add to the above Laplacian the
matrix $\lambda exp(\mathbf{U})$ and call the routine MATLAB \texttt{fsolve.}

\paragraph{Fist eigenfunction as an initial guess for Newton-Kantorovich algorithm.}
The first eigenfunction of the eigenproblem
\begin{equation}
\left\{ 
\begin{array}{c}
\Delta u+\lambda u=0,\ \textrm{in}\ \Omega :=\left( -L,\ L\right) \times\left( -L,\
L\right) , \\ 
u=0,\ \textrm{on}\ \partial \Omega .%
\end{array}%
\right.   \label{2DL_eq}
\end{equation}
will provide the initial guess for the following iterative process designed to solve the nonlinear algebraic system.

In order to solve the eigenproblem (\ref{2DL_eq}), i.e., to discretize the Laplacian operator, we use ChC. Subsequently we make use of MATLAB code \texttt{eig}. We find at least the first ten eigenvalues with a precision of order $O\left(10^{-10}\right).$
It is known that domains with symmetries will often have eigenvalues with multiplicities greater than one. The authors of \cite{KS} have observed that, for instance for the square of side $\pi,$ the eigenproblem (\ref{2DL_eq}) has
eigenvalues
$$ \lambda_{m,n}= m^2 + n^2,\ m,n=1,2,\ldots . $$
Whenever $m \neq n$, $\lambda_{m,n}$ will have multiplicity at least two. A simple shift of the domain $[-1, 1] \times [-1, 1]$ into $[0, \pi] \times [0, \pi]$ has enabled us to compare our computed eigenvalues with the above $ \lambda_{m,n}$. Thus we prove the accuracy of our computation.

\begin{figure*}
  \includegraphics[width=0.75\textwidth]{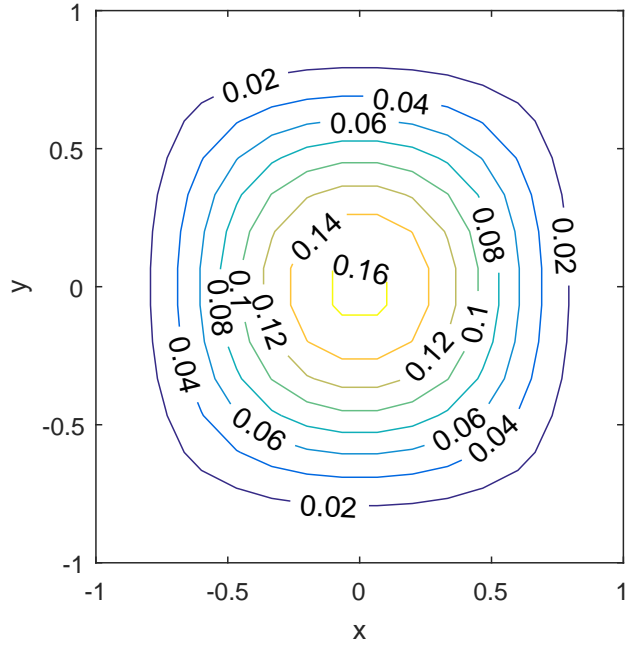}
\caption{Izolines of the \textquotedblleft{small}\textquotedblright solution to problem (\ref{2DB_eq}) produced by ChC with $N:=16$ when $\lambda:=0.5.$}
\label{fig:6}       
\end{figure*}

A contour plot which displays the isolines of the matrix solution to problem (\ref{2DB_eq}) is reported in Fig. \ref{fig:6}. It confirm the symmetries discussed above. 

\begin{figure}
  \includegraphics{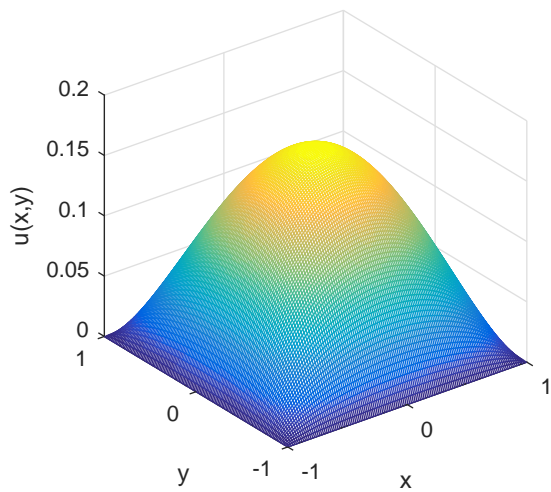}
\caption{ChC \textquotedblleft{small}\textquotedblright solution to 2D Bratu's problem for $\lambda=0.5$ when 
the order of approximation of ChC is $N:=20.$}
\label{fig:4}       
\end{figure}

In Fig. \ref{fig:4} we display a sixteen square-point approximation to problem (\ref{2DB_eq}). Actually, in order to have a more suggestive image we interpolate the solution to a finer grid and then plot. We use the MATLAB routine \texttt{interp2} with the option \texttt{'spline'}.

\begin{figure*}
  \includegraphics[width=0.75\textwidth]{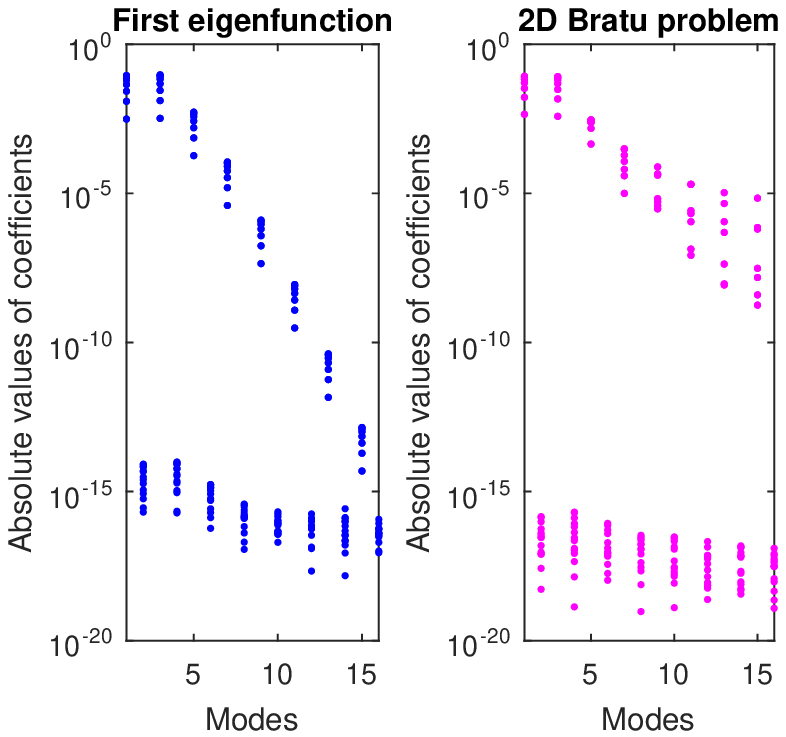}
\caption{Absolute values of the Chebyshev coefficients for the first eigenfunction (left panel) and the Chebyshev coefficients of solution to problem (\ref{2DB_eq}) (right panel). In both cases the order of approximation of ChC is $N:=16.$}
\label{fig:5}       
\end{figure*}

The Chebyshev coefficients of the first eigenfunction and of the solution to problem (\ref{1DB_eq}) are displayed in the left respectively right panels of Fig. \ref{fig:5} in a log-linear plot. They both confirm the assertion from Paragraph \ref{sing_sim}, i.e., the coefficients of the odd-degree Chebyshev polynomials must vanish. Actually, they are found in strips of points of order $10^{-15}$ or less. Comparing the two panels in this figure we also can conclude that the coefficients of the first eigenfunction fall off exponentially rather than algebraically with $N$ but the non linearity in Bratu's equation  reduces the rate of this fall. These coefficients are obtained using fast Chebyshev transform (FCT).

\paragraph{The lowest basis function as an initial guess for Newton-Kantorovich algorithm.}

Let's consider now the simplest basis function (\ref{1point_approx}) as an initial guess for Newton-Kantorovich algorithm. We want to approximate a solution situated on the right branch of the bifurcation diagram. 
\begin{figure}
  \includegraphics{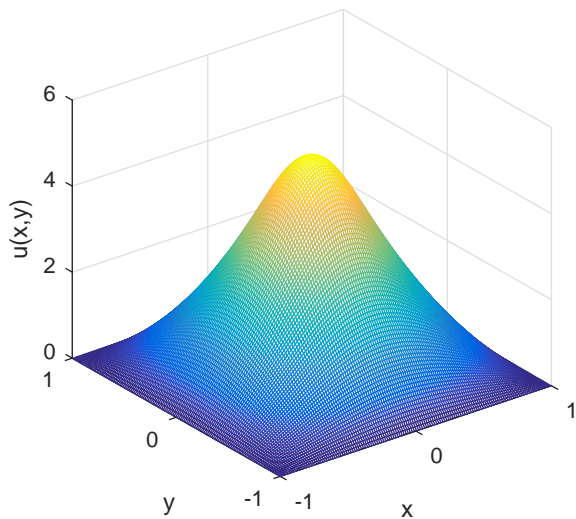}
\caption{ChC \textquotedblleft{big}\textquotedblright solution to 2D Bratu's problem for $\lambda=0.5$ when the order of approximation of ChC is $N:=16.$}
\label{fig:7}       
\end{figure}
Now the highest point of solution attains $$u_{max}= 4.677164395529806e+00$$ which is much higher than in the previous case when we have obtained only  $$u_{max}=1.865174060688610e-01.$$ It is also easy to see that the allure of the \textquotedblleft{big}\textquotedblright solution is much sharper than in the previous case.

\begin{figure*}
  \includegraphics[width=0.75\textwidth]{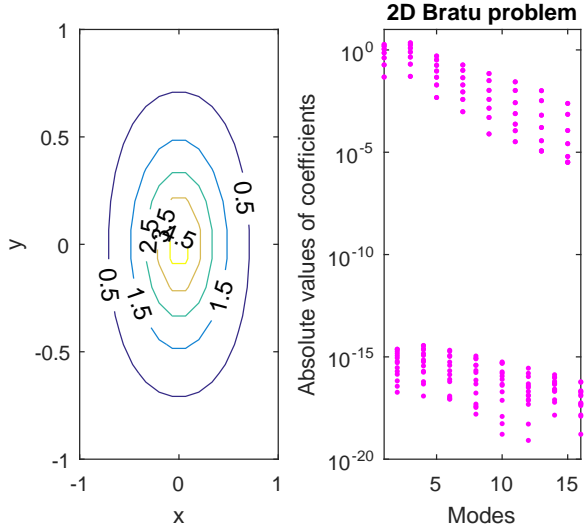}
\caption{Izolines of the \textquotedblleft{big}\textquotedblright solution to problem (\ref{2DB_eq}) (left panel) and in a log-linear plot the Chebyshev coefficients
of solution (right panel). All these results are produced by ChC when the order of approximation of ChC is $N:=16.$}
\label{fig:8}       
\end{figure*}

A contour plot which displays the isolines of the matrix solution in this case is reported in the left panel of Fig. \ref{fig:8}. It confirms the symmetries discussed above. 
In Table \ref{tab:1} we report some output parameters that \texttt{fsolve} provides when they solve the nonlinear algebraic system generated by ChC algorithm. Comparing its two lines we notice that the computation effort in the second case is double compared to those for the first case. In terms of elapsed time the difference is not so large.
\begin{table}
\caption{\texttt{fsolve} outcomes for both initial guesses.}
\label{tab:1}       
\begin{tabular}{llll}
\hline\noalign{\smallskip}
Initial guess & Nr. iter. & Func-count. & First-order optimality  \\
\noalign{\smallskip}\hline\noalign{\smallskip}
First eigenfunction of (\ref{2DL_eq}) & 4 & 985 & 2.73e-09\\
Lowest basis function (\ref{1point_approx})& 8 & 1773 & 2.05e-10\\
\noalign{\smallskip}\hline
\end{tabular}
\end{table}
\subsection{ChC solutions to some L-B-G problems.} \label{2DG-solution}
Let's consider in this section first the so called Gelfand's perturbed problem
 \begin{equation}
\left\{ 
\begin{array}{c}
\Delta u+\lambda e^{\left( \frac{u}{1+\varepsilon u}\right) }=0,\ \textrm{in}\ \Omega
:=\left( -L,\ L\right) \times \left( -L,\ L\right) ,\ 0<\varepsilon <<1, \\ 
u=0,\ \textrm{on}\ \partial \Omega .%
\end{array}%
\right.   \label{2DG_eq}
\end{equation}
We apply to this problem the same ChC technique along with Newton-Kantorovich and we do not notice major differences from the Bratu's problem regardless of the values of parameter $\varepsilon$.
The same conclusion holds when we consider in the equation (\ref{exothermic}) the non linearity $f(u):=\cosh(u)$ or $f(u):=\sinh(u).$ 
For the sake of brevity, we will not present results for these cases.
\section{Concluding remarks}

All in all, we have solved some 2D nonlinear Liouville-Bratu-Gelfand problems with exponential accuracy in
spite of corner singularities of the domain. The exponential accuracy is proved by the rate the Chebyshev coefficients of solutions decay. We have found out both, the \textquotedblleft{small}\textquotedblright solution as well as the \textquotedblleft{big}\textquotedblright solution. This latter result is rarely reported.

The Chebyshev pseudospectral method along with the Newton-Kantorovich algorithm is so efficient that it produce remarkably accurate results with a very modest approximation orders. Actually we have worked with $10<N<20$ and the iterative algorithm converged in just a few iterations. Also we have to remark that the Chebfun, used here for 1D problems, as well as MATLAB codes based on classical ChC are extremely simple, robust and efficient.

\end{document}